\documentclass[11pt]{amsart}

\usepackage{graphicx}
\usepackage{amsmath}
\usepackage{amssymb}
\usepackage{mathrsfs} 
\usepackage{xcolor}
\usepackage{bm}
\usepackage{hyperref}
\usepackage{booktabs}

\usepackage[utf8]{inputenc}

\newtheorem{thm}{Theorem}[section]
\newtheorem{prop}[thm]{Proposition}
\newtheorem{lem}[thm]{Lemma}

\theoremstyle{definition}
\newtheorem{definition}[thm]{Definition}

\theoremstyle{remark}

\newtheorem{remark}[thm]{Remark}

\newcommand{\R}{\mathbb{R}}  
\newcommand{\C}{\mathbb{C}} 
\newcommand{\Z}{\mathbb{Z}}

\newcommand{\F}{\mathbf{F}}

\DeclareMathOperator{\supp}{supp}

\newcommand{\ve}{\varepsilon}
\newcommand{\norm}[1]{\lVert #1 \rVert}
\newcommand{\abs}[1]{\lvert #1 \rvert}

\newcommand{\Fdft}{\mathbf{F}}

\begin{document}

\title[On the Donoho--Stark uncertainty principles]{On approximate support versions of Donoho--Stark uncertainty principles for $k$-Hadamard operators }

\subjclass[2020]{42A38, 42A99, 15B34, 94A12}
\keywords{Uncertainty principles, $k$-Hadamard matrices, Donoho--Stark uncertainty principle}

\author{Rami Ayoush}

\address{Faculty of Mathematics and Information Sciences, Warsaw University of Technology, ul.Koszykowa 75, 00-662 Warsaw, Poland}

\email{rami.ayoush@pw.edu.pl}

\address{Institute of Analysis, Johannes Kepler University Linz, Altenberger Strasse 69, A-4040 Linz, Austria}

\begin{abstract}
In this paper we present results related to Donoho--Stark type uncertainty principles that were proved in \cite{WW} by Wigderson and Wigderson for $k$-Hadamard operators in the finite-dimensional setting. Namely, the authors showed that Donoho-Stark principles can be generalized to versions for approximate supports with respect to norms $\ell^{p}$ for $p=1$ and $p=2$. In the second case they assumed the so-called speciality condition, namely that a suitable matrix satisfies $A^{*}A = kI$. First we show that such theorems cannot be generalized to $p>1$ if matrices are not special $k$-Hadamard. Next we prove Donoho-Stark type principle for special $k$-Hadamard operators for $\ell^{p}$ norms within the range $1\leq p\leq 2$. Finally, we show that such theorems are false for $p>2$ even in the presence of the speciality condition.
\end{abstract}

\maketitle


\section{Introduction} 

The theorem of Donoho and Stark (see \cite{DS}) says that for any function $f\colon \Z_{n} \to \C$ that is non-identically zero\footnote{Throughout the paper, functions and vectors are assumed to be non-zero.} we have 
\begin{equation}\label{DS-classical}
|\supp (f)| |\supp(\widehat{f})| \geq n,
\end{equation}
where the equality holds for $f$ being a specific `spike train', i.e. a function obtained from the indicator of a subgroup of $\Z_{n}$ by operations of shifting in the time and frequency domains, or scalar multiplication (see Appendix A in \cite{DS} for the details). 
The above estimate is also known to be true if we replace supports with certain other measures of concentration, as we describe below.

\begin{definition}\label{def:supp}
	For $f\colon \Z_{n} \to \C$, $T\subseteq \Z_n$, $p\in[1,\infty]$ and $\ve\in[0,1]$, we say that $f$ is
	\emph{$(p,\ve)$-supported on $T$} if $\norm{f\mathbf 1_{T^{c}}}_p\le\ve\norm f_p$, where $\mathbf 1_S$ denotes
	the indicator function of $S$ and $T^c=\Z_{n}\setminus T$. Put
	\[
	\abs{\supp^p_\ve(f)}=\min\bigl\{\abs T:\ f\text{ is }(p,\ve)\text{-supported on }T\bigr\}.
	\]
\end{definition}

\begin{definition}\label{def:band-lim}
	Let $p\in[1,\infty]$, $W \subseteq \Z_{n}$ and $Q=\mathcal F^{-1}\mathbf 1_W\mathcal F$ be the band-limiting operator. The function $f\colon \Z_{n} \to \C$ is said to be \emph{$(p,\delta)$-bandlimited to $W$} if there is a $g$ with
	$Qg=g$ and $\norm{f-g}_p\le\delta\norm f_p$.
\end{definition}

The notions of support, approximate support and band-limitedness may be also atributes of vectors $v \in \C^{n}$: one may simply consider them as functions on $[n] = \{0,...,n-1\}$.

The two classical Donoho--Stark approximate uncertainty principles on the cyclic group correspond to the exponents $p=1$ and $p=2$ in Definition~\ref{def:supp} and Definition~\ref{def:band-lim}.

\begin{thm}[\cite{DS}, Thm.~3]\label{thm:DS2}
	Let $f\colon \Z_{n} \to \C$, and let $\ve,\eta\in[0,1]$ with $\ve+\eta\leq1$. Then
	\[
	\abs{\supp^2_\ve(f)}\,\abs{\supp^2_\eta(\widehat f)}\ \ge\ n\,(1-\ve-\eta)^2 .
	\]
\end{thm}

The second one concerns $\ell_{1}$ norm and error measured on the time domain side. 

\begin{thm}[\cite{DS}, Thm.~7]\label{thm:DS1}
	Let $f\colon \Z_{n} \to \C$ and let $\ve,\eta\in[0,1)$. Suppose $f$ is
	$(1,\ve)$-supported on $T\subseteq\Z_n$ and $(1,\eta)$-bandlimited to
	$W\subseteq\Z_n$. Then
	\[
	\abs T\,\abs W\ \ge\ \frac{1-\ve-\eta}{1-\eta}\, n.
	\]
\end{thm}

In the literature there were several attempts to axiomatize various uncertainty principles (see e.g. \cite{EB}, \cite{GM}, \cite{GG} \cite{GJ}, \cite{GJ2}, \cite{JLW}, \cite{WW}, \cite{HLW}, \cite{DLP}), also outside the translation invariant setting. For instance, in the paper \cite{WW} whose results we study here, Wigderson and Wigderson isolated a simple criterion which guarantees that a broad family of uncertainty principles (including Donoho--Stark type principles) known already for the Fourier transform hold for any operator satisfying it. See `The primary uncertainty principle' in Section 2 in \cite{WW}. The main target of this methodology is the class of $k$-Hadamard operators whose definition we recall for the finite dimensional setting:

\begin{definition}\label{def:kH}
	Let $A\in\C^{m\times n}$ and $k>0$. Then $A$ is \emph{$k$-Hadamard} if it satisfies $\norm A_{1\to\infty}\le1$, the Gram matrix $G:=A^*A$ is invertible, and
	$\norm{G^{-1}}_{\infty\to\infty}\le1/k$. 
\end{definition}

There is an alternative formulation: $A$ is $k$-Hadamard if every entry of $A$ has modulus at most $1$ (equivalently $\norm A_{1\to\infty}\le1$), and
$\norm{A^*Av}_\infty\ge k\norm v_\infty$ for every $v\in\C^n$ (cf. Definition 2.2 in \cite{WW}).

\begin{definition}\label{def:specialkH}
	If $A$ is a square $k$-Hadamard matrix and
	$
	A^*A = kI,
	$
 	then we call $A$ \emph{special $k$-Hadamard} operator.
\end{definition}
The notion of $k$-Hadamard operators was introduced in \cite{WW}, while the speciality assumption is an adaptation of a similar notion from \cite{DLP} to the finite dimensional case. Examples of such matrices include DFT matrices, Walsh matrices, weighing matrices (always special $k$-Hadamard) and incidence matrices of finite geometries (non-special).  

By using their unified approach, the authors of \cite{WW} proved for this family, among other things, the following bounds. Their analog of \eqref{DS-classical}, Theorem~3.2 in \cite{WW}, states that
\begin{equation}
	|\supp (v)| |\supp(Av)| \geq k,
\end{equation} 
for every $k$-Hadamard $A$ and every $v\neq0$.

Another result, Theorem~3.17 therein, is an $\ell^1$-version for approximate supports which asserts that
\begin{equation}\label{eq:L1}
	\abs{\supp^1_\ve(v)}\,\abs{\supp^1_\eta(Av)}\ \ge\ k(1-\ve)(1-\eta) \ge\ k(1-\ve-\eta)
\end{equation}
for every $k$-Hadamard $A$.

Moreover, they proved a Euclidean version (Theorem~3.20 in \cite{WW}) that gives
\begin{equation}\label{eq:WW320}
	\abs{\supp^2_\ve(v)}\,\abs{\supp^2_\eta(Av)}\ \ge\ k(1-\ve-\eta)^2
	\qquad(\ve+\eta\le1)
\end{equation}
under the speciality assumption.

The additional speciality condition is necessary for the proof in \cite{WW} and the authors remarked that a potential programme of getting rid of it in $\eqref{eq:WW320}$ should require a method essentially different from theirs. It turns out, however, that such a general bound cannot be true due to the following result: 

\begin{thm}\label{thm:main}
	Let $1 < p \leq 2$.  Then, there is no function $\Phi:(0,\infty)\to(0,\infty)$ with $\Phi(k)\to\infty$ as $k\to\infty$
	such that
	\[
	\abs{\supp^p_\ve(v)}\,\abs{\supp^p_\eta(Av)}\ \ge\ \Phi(k)\,(1-\ve-\eta)^p
	\]
	holds for every $k>0$, every square $k$-Hadamard matrix $A$, every $v\neq0$ and all
	$\ve,\eta\in[0,1]$ with $\ve+\eta\le1$.
\end{thm}

In Section ~\ref{counterex} we present two constructions related to this theorem. The first one is a slightly simpler counterexample for $p>\frac{4}{3}$ which covers the Euclidean case. Its idea concerns building a sequence of operators that interpolate between discrete Fourier transforms and certain rank-one operators. The rank-one part produces the concentration, while the DFT part keeps the matrix $k$-Hadamard with $k \to +\infty$. In the next construction we prove the claim for other norms: we construct $k$-Hadamard operators from block-diagonal sums of DFT matrices. The aim of this operation is to obtain matrices with sparse rows which reinforces concentration behaviour in the image and allows one to overcome limitations of the previous example. 
	
Next we show a positive result, which is the following instance of the uncertainty principle: 
\begin{thm}\label{main_result}
Let $A\in\C^{n\times n}$ be a special $k$-Hadamard operator. Then for
every $p\in[1,2]$, every $v\in\C^n\setminus\{0\}$ and all $\ve,\eta\in[0,1]$ with
$\ve+\eta\le1$,
\begin{equation*}\label{eq:main}
	\abs{\supp^p_\ve(v)}\,\abs{\supp^p_\eta(Av)}\ \ge\ k\,(1-\ve-\eta)^p .
\end{equation*}	
\end{thm}

We divide the proof into two parts, both are contained in Section  ~\ref{specialUP}. The first involves a combination of methods derived from the papers \cite{DS} and \cite{Smith}, and covers the range $1\leq p \leq p_{B} \approx 1.85$. Here, the constant $p_B$ is a solution to the equation 
\[
p(2-p)-2^{-p} = 0, \quad p \in (1,2).
\] 
The second one deals with the range $\frac{5}{3}<p\leq 2$ and is based on a refinement of estimates from the Landau-Pollak article (\cite{LP}).

Finally, we prove the fact (which seems to be of a folklore flavour) that for all $p>2$ this theorem fails to hold in the approximate support form, even in the presence of the speciality criterion (we check it on a sequence of discrete Fourier transforms):

\begin{prop}\label{remark}
	Let $p>2$.  Then, there is no function $\Phi:(0,\infty)\to(0,\infty)$ with $\Phi(k)\to\infty$ as $k\to\infty$
	such that
	\[
	\abs{\supp^p_\ve(v)}\,\abs{\supp^p_\eta(Av)}\ \ge\ \Phi(k)\,(1-\ve-\eta)^p
	\]
	holds for every $k>0$, every special $k$-Hadamard matrix $A$, every $v\neq0$ and all
	$\ve,\eta\in[0,1]$ with $\ve+\eta\le1$.
\end{prop}

The above proposition is proved in Section ~\ref{folklore} in a slightly strengthened form.

We conclude the article with remarks explaining the relation of this research with \cite{Smith} and a summary of numerical experiments aimed at finding matrices related to Theorem ~\ref{thm:main} that are as small as possible.

\section{Notation and preliminaries}

Throughout the paper, we abbreviate $[n] = \{0,\dots,n-1\}$, and
$E\subseteq[n]$ and $\Sigma\subseteq[m]$ are any (possibly empty) sets. By $P_E$ and $P_\Sigma$ we denote the
corresponding coordinate orthogonal projections, $E^c$ is the complement of $E$, while $v[E]$ is the
restriction of $v$ to $E$. The symbol $\mathbf 1$ stands for the all-ones vector, $J=\mathbf1\mathbf1^{\!\top}$,
and $e_i$ is the $i$th standard basis vector. 

For $p\in[1,\infty]$ and $v\neq0$ we set
\begin{equation}\label{eq:params}
	\ve^{(p)}_E(v)=\frac{\norm{v-P_Ev}_p}{\norm v_p}=\frac{\norm{v[E^c]}_p}{\norm v_p},
	\qquad
	\eta^{(p)}_\Sigma(Av)=\frac{\norm{Av-P_\Sigma Av}_p}{\norm{Av}_p},
\end{equation}
where in the case $p=\infty$ the norms are just the usual suprema. The quantity $\eta^{(p)}_\Sigma$ is well defined because $k$-Hadamard operators are injective (which follows from Definition~\ref{def:kH} since $A^*A$ is invertible).

For every $p\in[1,\infty]$ and every $\ve\in[0,1]$ we denote
\begin{equation}\label{eq:suppchar}
	\abs{\supp^p_\ve(v)}=\min\bigl\{\abs E:\ \ve^{(p)}_E(v)\le\ve\bigr\},
\end{equation}
and likewise for $\abs{\supp^p_\eta(Av)}$ with $\eta^{(p)}_\Sigma$. When the sets $E,\Sigma$
are clear from the context we abbreviate
\[
\ve_p=\ve^{(p)}_E(v),\qquad \eta_p=\eta^{(p)}_\Sigma(Av).
\]

For $n\ge2$ let $\omega=e^{-2\pi i/n}$ and let $F_n$ be the unitary discrete Fourier transform
of $\Z_n$, $(F_n)_{jl}=n^{-1/2}\omega^{jl}$ for $j,l\in\Z_n$. Then
\begin{equation}\label{eq:dft}
	\Fdft_n:=\sqrt n\,F_n,\qquad (\Fdft_n)_{jl}=\omega^{jl},
\end{equation}
has all entries of modulus $1$ and satisfies $\Fdft_n^*\Fdft_n=\Fdft_n\Fdft_n^*=nI$, so
$\Fdft_n$ is a special $n$-Hadamard matrix.

In this article we use principal branches of $\arcsin$ and $\arccos$.
\section{Proof of Theorem \ref{thm:main}}\label{counterex}

\subsection{First construction: the interpolation trick}
Before we prove Theorem \ref{thm:main} in full generality we present a slightly simpler construction that works for exponents $p > \frac{4}{3}$.

For $\theta\in(0,1]$ let $D_\theta=\operatorname{diag}(1,\theta,\dots,\theta)$ and
\begin{equation}\label{eq:family}
	A_{n,\theta}:=D_\theta\F_n=\theta\F_n+(1-\theta)\,e_0\mathbf1^{\!\top},
\end{equation}
be a rank-one
perturbation of $\theta\F_n$ supported on one row. Here, $\mathbf 1$ is the all-ones column vector of size $n$. For $R>0$ put
\begin{equation}\label{eq:vec}
	v_R:=R\,e_0+\mathbf 1\in\C^n,\qquad E=\Sigma=\{0\},\quad\text{so } \abs E\abs\Sigma=1 .
\end{equation}

First we show that this construction yields $k$-Hadamard matrices.

\begin{lem}\label{lem:gram}
	Let $\theta\in(0,1]$ and $A=A_{n,\theta}$. Then $\norm A_{1\to\infty}\le1$ and
	\begin{equation}\label{eq:gram}
		A^*A=\theta^2nI+(1-\theta^2)J,\qquad J=\mathbf1\mathbf1^{\!\top},
	\end{equation}
	with eigenvalues $\theta^2n$ (multiplicity $n-1$) and $n$ (multiplicity $1$). Consequently
	$A$ is $k$-Hadamard for some $k$ satisfying
	\[
	k \ge \frac{\theta^2n}{2}.
	\]
	
\end{lem}

\begin{proof}
	Entries of $A$ are convex combinations of numbers of modulus at most $1$ and so their modulus is at
	most $1$ as well. Since $D_\theta^*D_\theta=\theta^2I+(1-\theta^2)e_0e_0^*$ and $e_0^*\F_n=\mathbf
	1^{\!\top}$, $\F_n^*e_0=\mathbf1$, we have
	\[
	A^*A=\F_n^*D_\theta^*D_\theta\F_n=\F_{n}^*(\theta^2\F_{n}+(1-\theta^2)e_{0}\mathbf1^{\!\top})
	=\theta^2nI+(1-\theta^2)J.
	\]
	Hence the spectrum of $A^*A$ is $\{\theta^2n,\ \theta^2n+(1-\theta^2)n\}=\{\theta^2n,n\}$ with the stated
	multiplicities. Indeed, $\mathbf{1}$ is an eigenvector corresponding to $n$ and $\mathbf1^{\!\perp}$ is the eigenspace related to $\theta^{2}n$. Since $J^2=nJ$, we can check that
	\[
	(A^*A)^{-1}=\frac{1}{\theta^2n}\Bigl(I-\frac{1-\theta^2}{n}J\Bigr),
	\]
	using $\theta^2n+(1-\theta^2)n=n$, and each row of $(A^*A)^{-1}$ has $\ell^1$ norm
	\[
	\frac{1}{\theta^2n}\Bigl(1-\frac{1-\theta^2}{n}+(n-1)\frac{1-\theta^2}{n}\Bigr)
	=\frac{1}{\theta^2n}\Bigl(1+(1-\theta^2)\frac{n-2}{n}\Bigr)
	\ \le\ \frac{2}{\theta^2n}.
	\]
	From the Definition~\ref{def:kH}, by inverting the last expression we get $k\ge\theta^2n/2$.
\end{proof}

\begin{remark}
	$A^*A=kI$ forces $(1-\theta^2)J=0$, i.e.
	$\theta=1$, hence for other parameters the matrix is not special $k$-Hadamard.
\end{remark}

Now let us determine suitable approximate supports.

\begin{lem}\label{lem:params}
	Let $A=A_{n,\theta}$, $v=v_R$ as in \eqref{eq:vec}. Then
	\begin{equation}\label{eq:image}
		Av=(R+n)\,e_0+\theta R\sum_{i\neq0}e_i ,
	\end{equation}
	and for every $p\in[1,\infty)$,
	
	\begin{equation}\label{eq:paramslem}
		\ve_p^p:=\bigl(\ve^{(p)}_E(v)\bigr)^p=\frac{n-1}{(R+1)^p+n-1},
	\end{equation}
	
	\begin{equation*}
		\eta_p^p:=\bigl(\eta^{(p)}_\Sigma(Av)\bigr)^p=\frac{(n-1)(\theta R)^p}{(R+n)^p+(n-1)(\theta R)^p}.
	\end{equation*}	

	Moreover $\abs{\supp^p_{\ve_p}(v)}=\abs{\supp^p_{\eta_p}(Av)}=1$.
\end{lem}

\begin{proof}
	The zeroth column of $\F_n$ is $\mathbf1$ and $\F_n\mathbf1=n e_0$, so
	$\F_nv=R\mathbf1+ne_0$; applying $D_\theta$ multiplies every coordinate but the zeroth by
	$\theta$, giving $(Av)_0=R+n$ and $(Av)_i=\theta R$ for $i\neq0$, which is \eqref{eq:image}.
	Now $v_0=R+1$ and $v_j=1$ for $j\neq0$, so $\norm v_p^p=(R+1)^p+(n-1)$ and
	$\norm{v[E^c]}_p^p=n-1$; likewise $\norm{Av}_p^p=(R+n)^p+(n-1)(\theta R)^p$ and
	$\norm{(Av)[\Sigma^c]}_p^p=(n-1)(\theta R)^p$. This is \eqref{eq:paramslem}. The support-size
	claims hold because $E,\Sigma$ are singletons and the defining inequalities hold with equality.
\end{proof}

We are now in a position to prove Theorem~\ref{thm:main} in the case $p>\frac{4}{3}$.

\begin{proof}[Proof of Theorem~\ref{thm:main} for $p>4/3$]
	Let $p>4/3$ and suppose such a $\Phi$ existed. Choose exponents $a,b$ with
	\begin{equation}\label{eq:exps}
		0<a<\tfrac12,\qquad \tfrac1p<b<a+1-\tfrac1p ,
	\end{equation}
	which is possible precisely because $\tfrac1p<a+1-\tfrac1p$ for some $a<\tfrac12$, i.e.\
	$\tfrac2p-1<\tfrac12$, i.e.\ $p>\tfrac43$. For $n\to\infty$ set $\theta=n^{-a}$, $R=n^{b}$ and
	$A=A_{n,\theta}$. By Lemma~\ref{lem:gram}, $A$ is $k$-Hadamard with
	$k\ge\theta^2n/2=\tfrac12n^{1-2a}\to\infty$. By
	\eqref{eq:paramslem}, we have
	\[
	\ve_p^p\lesssim \frac{n}{R^p}=n^{1-pb}\to0,
	\qquad
	\eta_p^p\lesssim \frac{n(\theta R)^p}{n^p}=n^{1-p+p(b-a)}\to0 ,
	\]
	where the convergence follows from the specific choice of parameters $a,b$. Hence $\ve_p+\eta_p\to0$, so
	$\ve_p+\eta_p\le1$ for large $n$, while by Lemma~\ref{lem:params} the left-hand side equals $1$
	and $\Phi(k)(1-\ve_p-\eta_p)^p\to\infty$: a contradiction.
\end{proof}

\subsection{Second construction: the block-diagonal counterexample}
In the previous construction, flatness of $\F_n$ seems to limit the first family to $p>4/3$. Replacing the
single DFT by a direct sum of small DFTs makes rows and columns sparse and removes this obstacle.

Let $N=sq$ with $s\ge2$, let $B=I_q\otimes\F_s\in\C^{N\times N}$ be block diagonal with $q$
copies of $\F_s$, and let $\iota\in\{0,1\}^N$ be the indicator of the set
$\{0,s,2s,\dots,(q-1)s\}$ of first indices of the blocks. Put
\begin{equation}\label{eq:block}
	C_{N,s}=\begin{pmatrix}\mathbf1^{\!\top}&-1\\[2pt] B&\iota\end{pmatrix}\in\C^{(N+1)\times(N+1)} ,
\end{equation}
and, for $R>0$, let $w_R\in\C^{N+1}$ have $(w_R)_0=R+1$, $(w_R)_j=1$ for $1\le j\le N-1$ and
$(w_R)_N=0$; take $E=\{0\}$ and $\Sigma=\{0\}$.

\begin{lem}\label{lem:block}
	With the above notation and $s^2\le N$:
	\begin{enumerate}
		\item[(i)] all entries of $C_{N,s}$ have modulus $0$ or $1$, and
		\[
		C_{N,s}^*C_{N,s}=\begin{pmatrix}sI_N+J_N&0\\0&1+N/s\end{pmatrix}.
		\]
		This notation means that in the last column and the last row all entries except the last one are zeros.
		\item[(ii)] $C_{N,s}$ is $k$-Hadamard with $k=s/2$.
		\item[(iii)] $(C_{N,s}w_R)_0=R+N$, and rows $1,\dots,N$ of $C_{N,s}w_R$ form the vector
		$R\,\mathbf1_{\text{block }0}+s\iota$; moreover
		$\abs{\supp^p_{\ve_p}(w_R)}=\abs{\supp^p_{\eta_p}(C_{N,s}w_R)}=1$ and
		\begin{equation}\label{eq:blockest}
			\ve_p := \ve^{(p)}_E(w_{R}) \le\frac{N^{1/p}}{R},\qquad
			\eta_p:= \eta^{(p)}_\Sigma(C_{N,s}w_{R}) \le\frac{R\,s^{1/p}}{N}+\Bigl(\frac sN\Bigr)^{1-\frac1p} .
		\end{equation}
	\end{enumerate}
\end{lem}

\begin{proof}
	(i) Writing $c=(-1,\iota)$ for the last column and
	$c_j=(1,Be_j)$ for $j<N$ (we enumerate columns from $0$ to $N$): for $i,j<N$, $\langle c_i,c_j\rangle=1+(B^*B)_{ij}=1+s\delta_{ij}$,
	since $B^*B=sI$; also $\norm c_2^2=1+q=1+N/s$. Finally, if $j$ lies in some block then $Be_j$ is
	supported in the same block and $\iota^{\!\top}Be_j$ is its entry at the block's first index, i.e.\
	the $j'$th entry, $j' \equiv j \mod s$ , of the first row of $\F_s$, which is $1$; hence
	$\langle c,c_j\rangle=-1+1=0$.
	
	(ii) By (i) the Gram matrix is block diagonal, so $\norm{(C^*C)^{-1}}_{\infty\to\infty}$ is the
	larger of the two blocks' values. The second is $(1+N/s)^{-1}\le s/N\le1/s$ by $s^2\le N$. For
	the first, $(sI_N+J_N)^{-1}=s^{-1}(I-\frac{1}{s+N}J)$ has rows of $\ell^1$ norm
	$s^{-1}(1+\frac{N-2}{s+N})\le2/s$. Hence $\norm{(C^*C)^{-1}}_{\infty\to\infty}\le2/s=1/k$.
	
	(iii) Write $w_R'=Re_0+\mathbf1_N$ for the vector consisting of first $N$ coordinates of $w_R$. Rows
	$1,\dots,N$ of $Cw_R$ are $Bw_R'+\iota\,(w_R)_N=Bw_R'$: the last column of $C$ contributes
	nothing, since $(w_R)_N=0$. Since the first column of $\F_s$ is $\mathbf1$ and
	$\F_s\mathbf1=s\,e_0$ within each block, $Be_0=\mathbf1_{\text{block }0}$ and
	$B\mathbf1_N=s\iota$. Hence
	$Bw_R'=R\,\mathbf1_{\text{block }0}+s\iota$. The adjoined row gives
	$(Cw_R)_0=\mathbf1^{\!\top}w_R'-1\cdot0=R+N$. Again, the support-size assertions hold because $E$ and
	$\Sigma$ are singletons and the defining inequalities hold with equality.
	
	For \eqref{eq:blockest}: $w_R[E^c]$ consists of $N-1$ entries equal to $1$ and one zero, so
	$\norm{w_R[E^c]}_p^p=N-1$, while $\norm{w_R}_p\ge\abs{(w_R)_0}=R+1$; hence
	$\ve_p^p\le(N-1)/(R+1)^p\le N/R^p$. As for $\eta_p$, the triangle inequality gives
	$\norm{(Cw_R)[\Sigma^c]}_p=\norm{R\,\mathbf1_{\text{block }0}+s\iota}_p
	\le R\norm{\mathbf1_{\text{block }0}}_p+s\norm\iota_p=Rs^{1/p}+s(N/s)^{1/p}$, on the other hand $\norm{Cw_R}_p\ge\abs{(Cw_R)_0}=R+N\ge N$; dividing those two quantities
	gives the second estimate.
\end{proof}
\begin{proof}[Proof of Theorem~\ref{thm:main}]
	Let $p>1$ and suppose such a $\Phi$ existed. Put
	\[
	\gamma=\tfrac12\min\{\tfrac12,\,p-1\}\in(0,\tfrac14],\qquad \tau=1-\frac{1+\gamma}{p}>0 ,
	\]
	where the positivity of $\tau$ follows from $\gamma\le(p-1)/2<p-1$. For $N\to\infty$ along powers of $2$ set
	$s=2^{\lceil\gamma\log_2N\rceil}$, $\delta=N^{-\tau/2}$ and $R=\delta Ns^{-1/p}$; for large $N$
	one has $2\le s\le\sqrt N$, so Lemma~\ref{lem:block} applies to $C=C_{N,s}$, which is
	$k$-Hadamard with $k=s/2\ge N^\gamma/2\to\infty$. By \eqref{eq:blockest} from the previous lemma and $s\le2N^\gamma$,
	\[
	\ve_p\le\frac{N^{1/p-1}s^{1/p}}{\delta}\le\frac{2^{1/p}N^{\frac{1+\gamma}{p}-1}}{\delta}
	=2^{1/p}\,\frac{N^{-\tau}}{\delta}=2^{1/p}\delta\longrightarrow0,
	\]
	\[
	\eta_p\le\delta+\Bigl(\frac sN\Bigr)^{1-\frac1p}\le\delta+\bigl(2N^{\gamma-1}\bigr)^{1-\frac1p}
	\longrightarrow0,
	\]
	since $\gamma<1$ and $1-1/p>0$. Hence $\ve_p+\eta_p\le1$ for large $N$, the left-hand side
	equals $1$ by Lemma~\ref{lem:block}(iii), and $\Phi(k)(1-\ve_p-\eta_p)^p\to\infty$:
	contradiction.
\end{proof}

\section{Proof of Theorem \ref{main_result}}\label{specialUP}

The main idea of the proof is to use two estimates, one near $p=1$ and the other one near $p=2$ that are sharper than Theorem \ref{main_result} for those parameters (see Remark \ref{sharpen1} and Remark \ref{sharpen2}), and try to apply them as far from those exponents as possible. In this way, we will prove the desired estimates on two overlapping intervals that cover the segment $[1,2]$.

Let $p\in[1,2]$ and $p'$ be the conjugate exponent ($p'=\infty$ when $p=1$), and
\[
\theta=\frac1p-\frac1{p'}=\frac2p-1\in[0,1].
\]
For $S\subseteq[n]$ let $P_S$ denote the coordinate projection. For a vector $x$ supported on
$S$ and $q\ge p$, Hölder's inequality gives
\begin{equation}\label{eq:holder}
	\norm{x}_p\le\abs S^{\frac1p-\frac1q}\norm x_q .
\end{equation}

\subsection{Range $1\leq p \leq p_{B} \approx 1.85$}\label{first_part}
Next Hausdorff-Young type estimate for $k$-Hadamard operators is essentially Proposition 3.24 in \cite{WW}. We reprove it here for completeness.
\begin{lem}\label{lem:HY}
	Let $A$ be special and $1 \leq p \leq 2$. Then for
	every $v$,
	\[
	\norm{Av}_{p'}\le k^{1/p'}\norm v_p ,
	\]
	and likewise for $A^*$ and $A^{\!\top}$ in place of $A$.
\end{lem}

\begin{proof}
	This is a standard interpolation argument. We have $\norm A_{1\to\infty}\le1$ and
	$\norm A_{2\to2}=\sqrt k$. The Riesz--Thorin theorem with $\frac1p=\frac{1-t}{1}+\frac t2$,
	$\frac1{p'}=\frac t2$, i.e.\ $t=2/p'$, gives $\norm A_{p\to p'}\le(\sqrt k)^{2/p'}=k^{1/p'}$.
	
	The same argument applies to $A^*$ and $A^{\!\top}$.
\end{proof}

For $\Sigma\subseteq[n]$ let us define an analog of band-limiting operator for our setting:
\begin{equation}\label{eq:Q}
	Q_\Sigma=k^{-1}A^*P_\Sigma A .
\end{equation}
By the speciality assumption and Lemma~\ref{lem:HY}, $Q_\Sigma$ is the orthogonal projection onto
$\{g:\ \supp(Ag)\subseteq\Sigma\}$. The following is Smith's estimate (Thm.~3.1 in \cite{Smith}) in the
present setting. Smith's proof uses the translation invariance of the Fourier transform on a
group, which is unavailable here; we replace it with a row-by-row bound.

\begin{lem}\label{lem:smith}
	Let $A$ be special, $1\leq p \leq 2$ and $E,\Sigma\subseteq[n]$. Then for every $v$ we have
	\[
	\norm{P_EQ_\Sigma v}_p\le\Bigl(\frac{\abs E\abs\Sigma}{k}\Bigr)^{1/p}\norm v_p .
	\]
\end{lem}

\begin{proof}
	Write $(P_EQ_\Sigma v)(x)=\sum_yK(x,y)v_y$. By \eqref{eq:Q},
	\[
	K(x,y)=\mathbf 1_E(x)\,k^{-1}\sum_{j\in\Sigma}\overline{a_{jx}}\,a_{jy},
	\]
	so for $x\in E$ the
	row $K(x,\cdot)$ equals $k^{-1}A^{\!\top}w_x$, where $(w_x)_j=\mathbf1_\Sigma(j)\overline{a_{jx}}$.
	Since $\abs{a_{jx}}\le1$, $\norm{w_x}_p\le\abs\Sigma^{1/p}$, and Lemma~\ref{lem:HY} for
	$A^{\!\top}$ gives $\norm{K(x,\cdot)}_{p'}\le k^{-1}k^{1/p'}\abs\Sigma^{1/p}
	=k^{-1/p}\abs\Sigma^{1/p}$. By Hölder,
	\[\abs{(P_EQ_\Sigma v)(x)}\le\mathbf1_E(x)\,k^{-1/p}\abs\Sigma^{1/p}\norm v_p; 
	\]
	raising to the
	power $p$ and summing over $x$ gives the claim.
\end{proof}

\begin{prop}\label{prop:M3}
	Let $A$ be special, $1\le p\le2$, $v\neq0$, and let $E,\Sigma\subseteq[n]$ satisfy
	\[
	\norm{v[E^c]}_p\le\ve\norm v_p,\qquad \norm{(Av)[\Sigma^c]}_p\le\eta\norm{Av}_p,
	\qquad 0\le\eta<1 .
	\]
	Then, with $X=\abs E\abs\Sigma/k$,
	\begin{equation}\label{eq:M3}
		X^{1/p}+\frac{\eta}{(1-\eta^p)^{1/p}}\,X^{\theta}\ \ge\ 1-\ve .
	\end{equation}
	The same holds with the roles of $\ve$ and $\eta$ interchanged, provided $\ve<1$:
	\begin{equation}\label{eq:M3p}
		X^{1/p}+\frac{\ve}{(1-\ve^p)^{1/p}}\,X^{\theta}\ \ge\ 1-\eta .
	\end{equation}
\end{prop}

\begin{proof}
	Since $k^{-1}A^*A=I$, we have $v-Q_\Sigma v=k^{-1}A^*P_{\Sigma^c}Av$, and hence
	\begin{equation}\label{decomposition}
	P_EQ_\Sigma v=v-P_{E^c}v-P_E\bigl(v-Q_\Sigma v\bigr).
	\end{equation}
	The middle term has norm at most $\ve\norm v_p$. For the last term, \eqref{eq:holder} on $E$
	with $q=p'$ and Lemma~\ref{lem:HY} for $A^*$ give
	\begin{multline}\label{third_term}
	\norm{P_E(v-Q_\Sigma v)}_p\le\abs E^\theta\norm{v-Q_\Sigma v}_{p'}
	\le\abs E^\theta k^{-1}k^{1/p'}\norm{(Av)[\Sigma^c]}_p \\
	\le\abs E^\theta k^{-1/p}\,\eta\norm{Av}_p .
	\end{multline}
	To return to $\norm v_p$, use the concentration of $Av$ on $\Sigma$: from
	\begin{multline}
	\norm{Av}_p^p=\norm{(Av)[\Sigma]}_p^p+\norm{(Av)[\Sigma^c]}_p^p\le
	\norm{(Av)[\Sigma]}_p^p+\eta^p\norm{Av}_p^p,
	\end{multline}
	together with \eqref{eq:holder} on $\Sigma$ and
	Lemma~\ref{lem:HY},
	\begin{multline}\label{norm_Av}
	\norm{Av}_p\le(1-\eta^p)^{-1/p}\norm{(Av)[\Sigma]}_p
	\le(1-\eta^p)^{-1/p}\abs\Sigma^\theta\norm{Av}_{p'} \\
	\le(1-\eta^p)^{-1/p}\abs\Sigma^\theta k^{1/p'}\norm v_p .
	\end{multline}
	By combining inequalities \eqref{third_term} and \eqref{norm_Av}, and using $k^{1/p'-1/p}=k^{-\theta}$,
	\[
	\norm{P_E(v-Q_\Sigma v)}_p\le\frac{\eta}{(1-\eta^p)^{1/p}}\,X^\theta\norm v_p .
	\]
	Therefore, by the above and \eqref{decomposition} we get 
	\[
	\norm{P_EQ_\Sigma v}_p\ge\bigl(1-\ve-\frac{\eta}{(1-\eta^p)^{1/p}}X^\theta\bigr)
	\norm v_p,
	\] 
	while Lemma~\ref{lem:smith} gives $\norm{P_EQ_\Sigma v}_p\le X^{1/p}\norm v_p$. Last two estimates prove
	\eqref{eq:M3}. 
	
	For \eqref{eq:M3p}, apply \eqref{eq:M3} to the special matrix $A^*$, the vector
	$u=Av$ and the sets $\Sigma,E$ in place of $E,\Sigma$: since $A^*u=kv$ and concentration is
	invariant under scaling, $A^*u$ is concentrated on $E$ with parameter $\ve$, and
	$\abs\Sigma\abs E/k=X$.
\end{proof}

\begin{remark}\label{sharpen1}
At $p=1$ we have $\theta=1$, and \eqref{eq:M3} reads $X\bigl(1+\frac{\eta}{1-\eta}\bigr)\ge1-\ve$,
i.e.\ $\abs E\abs\Sigma\ge k(1-\ve)(1-\eta)$ which is the left inequality in \eqref{eq:L1}.
\end{remark}

Next lemma provides a range of applicability of Proposition ~\ref{prop:M3} to our problem. It turns out that it effectively works for exponents starting from $p=1$ up to the exponent that it is determined by the root of the equation
\[
p(2-p)-2^{-p} = 0, \quad p \in (1,2).
\]
We denote this root by $p_B$. Numerical check shows that $p_B \approx 1.85$.

\begin{lem}\label{lem:elem}
	Let $1\le p\le p_B$ and $\sigma=p(2-p)$. Then for every $m\in[0,\frac12]$,
	\[
	(1-2m)^{\sigma}\ \le\ 1-m^p .
	\]
\end{lem}

\begin{proof}
	Put $h(p)=p(2-p)-2^{-p}$. Then $h(1)=\frac12$, $h(2)=-\frac14$ and
	$h''(p)=-2-2^{-p}(\log2)^2<0$, so $h$ is concave with a unique root $p_B\in(1,2)$, and $h\ge0$
	on $[1,p_B]$, i.e.\ $2^{-p}\le\sigma$ there. Since $0<\sigma\le1$, Bernoulli inequality gives $(1-x)^\sigma\le1-\sigma x$ for $x\in[0,1]$; with $x=2m$,
	$(1-2m)^\sigma\le1-2\sigma m$. Finally $m^{p-1}\le2^{1-p}\le2\sigma$ for $m\le\frac12$, so
	$m^p\le2\sigma m$ and $1-2\sigma m\le1-m^p$.
\end{proof}

\begin{proof}[Proof of Theorem~\ref{main_result} for $1\leq p \leq p_B$]
	The case $\ve+\eta=1$ is trivial, so assume $\ve+\eta<1$. Let us denote
	$s=1-\ve-\eta\in(0,1]$, $\mu=\min\{\ve,\eta\}\in[0,\frac12)$ and $M=\max\{\ve,\eta\}$. Let $E$
	and $\Sigma$ realize $\abs{\supp^p_\ve(v)}$ and $\abs{\supp^p_\eta(Av)}$, and put
	$X=\abs E\abs\Sigma/k$. Applying \eqref{eq:M3} if $\eta\le\ve$ and \eqref{eq:M3p} otherwise, we get
	\begin{equation}\label{eq:sym}
		f(X)\ge1-M,
	\end{equation}
	where
	\begin{equation*}
	f(Y):=Y^{1/p}+c\,Y^\theta,\qquad c=\frac{\mu}{(1-\mu^p)^{1/p}}.
	\end{equation*}
	We will argue by contradiction. Suppose that $X<s^p$. Since $f$ is strictly increasing on $[0,\infty)$,
	\[1-M\le f(X)<f(s^p)=s+c\,s^{p\theta}=s+c\,s^{2-p}.\] 
	As $1-M-s=\mu$, this says
	$\mu<c\,s^{2-p}$. If $\mu=0$ then $c=0$, which is absurd. Otherwise it says
	$(1-\mu^p)^{1/p}<s^{2-p}$, i.e.\ $s^{\sigma}>1-\mu^p$ with $\sigma=p(2-p)$. But
	$s=1-\mu-M\le1-2\mu$, so $(1-2\mu)^\sigma\ge s^\sigma>1-\mu^p$, contradicting
	Lemma~\ref{lem:elem}. Hence $\abs E\abs\Sigma\ge ks^p$.
\end{proof}

\subsection{Range $\frac{5}{3}< p \leq 2$}\label{second_part}
Below we start with a lemma that is useful for comparing concentration with respect to various norms. More specifically, the estimates in this subsection will be based on bounds that are most effective at $p=2$ and the role of this lemma will be to incorporate them to smaller exponents.

\begin{lem}\label{lem:comp}
	Let $v\neq0$, $1\le p\le2$ and $\ve\in[0,1]$. Then
	$\abs{\supp^2_{\ve^{p/2}}(v)}\le\abs{\supp^p_\ve(v)}$.
\end{lem}

\begin{proof}
	Let $m=\abs{\supp^p_\ve(v)}$ and let $T$ be a set of $m$ indices carrying the $m$ largest values
	of $\abs{v_j}$. Among sets of size $m$, $T$ minimizes $\norm{v[T^c]}_s$ for every $s$, so $v$ is
	$(p,\ve)$-supported on $T$. For $s>0$ put $S_s=\sum_{j\in T}\abs{v_j}^s$ and
	$R_s=\sum_{j\notin T}\abs{v_j}^s$, and let $\mu=\min_{j\in T}\abs{v_j}\ge
	\nu=\max_{j\notin T}\abs{v_j}$. Suppose first that $R_2\neq0$. We have $S_2\ge
	\mu^{2-p}S_p$ and $R_2\le\nu^{2-p}R_p$, so $S_2/R_2\ge(\mu/\nu)^{2-p}S_p/R_p\ge S_p/R_p$, and
	\[
	\frac{\norm{v[T^c]}_2^2}{\norm v_2^2}=\frac{R_2}{S_2+R_2}\le\frac{R_p}{S_p+R_p}
	=\frac{\norm{v[T^c]}_p^p}{\norm v_p^p}\le\ve^p .
	\]
	The above inequality is trivial when $R_2=0$. Hence $v$ is
	$(2,\varepsilon^{p/2})$-supported on $T$, and $|T|=m$ gives the claim.
\end{proof}

The next lemma is essentially a part of Theorem 2 in \cite{LP} stated here in a suitable form. Below we use principal branches of $\arcsin$ and $\arccos$.

\begin{lem}\label{lem:LP}
	Let $A$ be special, $v\neq0$, and $a,b\in[0,1]$ with $\arcsin a+\arcsin b\le\pi/2$. Then
	\[
	\abs{\supp^2_a(v)}\,\abs{\supp^2_b(Av)}\ \ge\ k\,\Theta(a,b)^2,\qquad
	\Theta(a,b)=\cos(\arcsin a+\arcsin b).
	\]
\end{lem}

\begin{proof}
	If $\Theta(a,b)=0$ there is nothing to prove, so assume $\arcsin a+\arcsin b<\pi/2$. Let $E$ and
	$\Sigma$ realise the two support sizes, put $U=k^{-1/2}A$, which is unitary, and normalise
	$\norm v_2=1$, so that $w=Uv$ is a unit vector which is $(2,b)$-supported on $\Sigma$. Let
	$\varphi_1=\arccos\norm{P_Ev}_2$ and $\varphi_2=\arccos\norm{P_\Sigma w}_2$; then
	$\sin\varphi_1\le a$ and $\sin\varphi_2\le b$, so $\varphi_1+\varphi_2<\pi/2$, and in particular
	$P_Ev\ne0\ne P_\Sigma w$. Put $z=P_Ev/\norm{P_Ev}_2$, $x=Uz$ and $y=P_\Sigma w/\norm{P_\Sigma
		w}_2$. On the unit sphere of $\C^n=\R^{2n}$ the geodesic distance is
	$d(\xi,\zeta)=\arccos\operatorname{Re}\langle\xi,\zeta\rangle$ (see \cite{BH}, Chapter I.2, Prop. 2.1. and Prop. 2.3.); since $U$ is unitary,
	$d(w,x)=d(v,z)=\varphi_1$, and $d(w,y)=\varphi_2$. The triangle inequality gives
	$d(x,y)\le\varphi_1+\varphi_2$, hence
	\[
	\norm{P_\Sigma UP_E}\ \ge\ \abs{\langle y,P_\Sigma UP_Ez\rangle}=\abs{\langle y,x\rangle}
	\ \ge\ \cos(\varphi_1+\varphi_2)\ \ge\ \Theta(a,b).
	\]
	On the other hand $\norm{P_\Sigma UP_E}\le\norm{P_\Sigma UP_E}_{HS}
	=\bigl(\sum_{i\in\Sigma,\,j\in E}\abs{a_{ij}}^2/k\bigr)^{1/2}\le(\abs E\abs\Sigma/k)^{1/2}$.
\end{proof}

\begin{remark}\label{sharpen2}
Let us observe that
\[
\Theta(a,b)=\sqrt{1-a^2}\sqrt{1-b^2}-ab.
\]
One can check that $\Theta(a,b)\ge1-a-b$, and therefore Lemma~\ref{lem:LP} sharpens \eqref{eq:WW320}.
\end{remark}

\begin{lem}\label{lem:q} For $r\ge1$ and $x\in[0,\pi/3]$ let us put
	\[
	q_r(x)=2^{1-r}\Bigl(1-\frac{x^2}{6}\Bigr)^r-\frac r2\,x^{2-r}.
	\]
	If $1\le r\le6/5$ and $0\le x\le\pi/3$, then $q_r(x)\ge0.059$.
\end{lem}

\begin{proof}
	The first term of $q_r$ is decreasing in $x$ and the second increasing, so $q_r(x)\ge
	q_r(\pi/3)$. The first term equals $2\bigl(\frac12(1-\frac{x^2}6)\bigr)^r$, which is decreasing
	in $r$; the second has $r$-derivative $\frac12x^{2-r}(1-r\log x)\ge0$, because
	$r\log x\le\frac65\log\frac\pi3<1$. Hence $q_r(x)\ge q_{6/5}(\pi/3)$. Finally,
	$1-\pi^2/54\ge0.8172$, $0.8172^{6/5}\ge0.784$, $2^{-1/5}\ge0.870$ and
	$\frac35(\pi/3)^{4/5}\le0.6228$, so $q_{6/5}(\pi/3)\ge0.784\cdot0.870-0.6228\ge0.059$.
\end{proof}

Further we will use three elementary numerical inequalities:

\begin{lem}\label{numerical}
	We have the following estimates:
		 \[\sin w\ge w(1-w^2/6)\ge0 \quad \text{for } 0\le w\le\pi/3,\]
		 \[\cos^ru\ge1-r(1-\cos u)\ge1-ru^2/2 \text{ for } 0\le u\le\pi/2 \text{ and } r\ge1,\]
		 \[s^r+t^r\ge2^{1-r}(s+t)^r \text{ for } s,t\ge0 \text{ and } r\ge1.\]
\end{lem}

\begin{proof}
	To get the first one it is enough to compare first derivatives, second follows from Bernoulli inequality, and the third one is an obvious consequence of the power means inequality.
\end{proof}

\begin{prop}\label{prop:trig}
	Let $1\le r\le6/5$ and $\alpha,\beta\ge0$ with $\alpha+\beta\le\pi/2$. Then
	\[
	\sin^r\alpha+\sin^r\beta+\cos^r(\alpha+\beta)\ \ge\ 1.
	\]
\end{prop}

\begin{proof}
	Write $S$ for the left-hand side and $\gamma=\alpha+\beta$.
	
	\emph{Case 1: $\gamma\le\pi/3$.} Since $\alpha,\beta\le\gamma$,
	$\sin^r\alpha+\sin^r\beta\ge(1-\gamma^2/6)^r(\alpha^r+\beta^r)\ge(1-\gamma^2/6)^r2^{1-r}\gamma^r$.
	Together with $\cos^r\gamma\ge1-r\gamma^2/2$ this gives $S-1\ge\gamma^r\,q_r(\gamma)\ge0$.
	
	\emph{Case 2: $\gamma>\pi/3$.} Then $\max\{\alpha,\beta\}>\pi/6$; by symmetry let $\alpha>\pi/6$.
	Put $u=\pi/2-\alpha\in[0,\pi/3)$ and $w=\beta$, so that $0\le w\le u$. Then
	$\sin\alpha=\cos u$ and $\cos(\alpha+\beta)=\sin(u-w)$, whence, using $w,u-w\le u$,
	\begin{multline}
	S=\cos^ru+\sin^rw+\sin^r(u-w)\ \ge\ 1-\frac r2u^2+\Bigl(1-\frac{u^2}6\Bigr)^r2^{1-r}u^r \\
	=1+u^r\,q_r(u)\ \ge\ 1,
	\end{multline}
	which proves the claim
\end{proof}

\begin{proof}[Proof of Theorem~\ref{main_result} for $p>\frac{5}{3}$]
	Let $r=2/p\in[1,6/5]$. If $\ve+\eta=1$ the right-hand side in Theorem ~\ref{main_result} vanishes, so let
	$\ve+\eta<1$. Put $a=\ve^{p/2}$, $b=\eta^{p/2}$, $\alpha=\arcsin a$, $\beta=\arcsin b$, so that
	$\ve=\sin^r\alpha$ and $\eta=\sin^r\beta$. Since $r\le2$,
	$\sin^2\alpha+\sin^2\beta\le\sin^r\alpha+\sin^r\beta=\ve+\eta\le1$, so $\sin^2\beta\le
	\cos^2\alpha = \sin(\frac{\pi}{2}-\alpha)$, i.e.\ $\alpha+\beta\le\pi/2$ (by the monotonicity of the sine on $[0, \pi/2]$). Hence $\Theta(a,b)=\cos(\alpha+\beta)\ge0$. Proposition~\ref{prop:trig} gives $\cos^r(\alpha+\beta)\ge1-\ve-\eta\ge0$, i.e.\
	$\Theta(a,b)\ge(1-\ve-\eta)^{p/2}$. By Lemma~\ref{lem:comp} applied to $v$ and to $Av$, and
	Lemma~\ref{lem:LP},
	\begin{multline}
	\abs{\supp^p_\ve(v)}\,\abs{\supp^p_\eta(Av)}\ \ge\ \abs{\supp^2_a(v)}\,\abs{\supp^2_b(Av)}
	\ \ge\ k\,\Theta(a,b)^2\ \\ \ge\ k\,(1-\ve-\eta)^p .
	\end{multline}
\end{proof}

\begin{proof}[Proof of Theorem ~\ref{main_result}]
	Intervals $[1,p_B]$ and $(5/3, 2]$ cover $[1,2]$, so by this and the previous subsection the theorem holds true for all parameters.
\end{proof}

\section{Proposition \ref{remark}: The case  $p>2$}\label{folklore}

In this section we prove that the approximate support Donoho--Stark theorem cannot hold for $p>2$ even for the DFT matrices.

The next lemma follows from direct computation.

\begin{lem}\label{lem:dftbound}
	With the notation of Lemma~\ref{lem:params} for $\theta = 1$ and the choice $R=\sqrt n$, for every $p>2$ we have
	\[
	\ve_p\leq n^{\frac1p-\frac12},\qquad \eta_p\leq n^{\frac1p-\frac12},
	\]
	and both quantities tend to zero as $n \to +\infty$.
\end{lem}

\begin{thm}\label{thm:pgt2}
	Let $p>2$. There is no function $\Phi:(0,\infty)\to(0,\infty)$ with $\Phi(k)\to\infty$ as
	$k\to\infty$ such that
	\[
	\abs{\supp^p_\ve(v)}\,\abs{\supp^p_\eta(Av)}\ \ge\ \Phi(k)\,(1-\ve-\eta)^p
	\]
	holds for every $k>0$, every special $k$-Hadamard matrix $A$, every
	$v\neq0$ and all $\ve,\eta\in[0,1]$ with $\ve+\eta\le1$. The same is true with
	$(1-\ve-\eta)^p$ replaced by any $\psi(\ve,\eta)$ with
	$\liminf_{(\ve,\eta)\to(0,0)}\psi(\ve,\eta)>0$.
\end{thm}

\begin{proof}
	Suppose such a $\Phi$ existed. For $n\ge2$ let $A=\Fdft_n$, which is $n$-Hadamard with
	$A^*A=nI$ by \eqref{eq:dft}, so $k=n$; let $v=\sqrt n\,e_0+\mathbf1$ and $E=\Sigma=\{0\}$. By
	Lemma~\ref{lem:dftbound},
	\[
	\ve_p+\eta_p\ \leq\ 2n^{\frac1p-\frac12}\ \longrightarrow\ 0\qquad(n\to\infty),
	\]
	since $p>2$; in particular $\ve_p+\eta_p\le1$ for $n$ large. By Lemma~\ref{lem:params} for $\theta = 1$ the
	left-hand side equals $1$ for every $n$, whereas
	\[
	\Phi(n)\bigl(1-\ve_p-\eta_p\bigr)^p\ \ge\ \Phi(n)\bigl(1-2n^{\frac1p-\frac12}\bigr)^p
	\ \longrightarrow\ \infty ,
	\]
	a contradiction for large $n$. The last assertion follows in the same way, since
	$\psi(\ve_p,\eta_p)$ is bounded below by a positive constant for $n$ large.
\end{proof}

\section{Further remarks}
\subsection{Relations to Smith's theorem}\label{sec:smith}

In the article \cite{Smith} Smith extended Theorem ~\ref{thm:DS2} and Theorem ~\ref{thm:DS1} to $L^p$,
$1\le p\le2$, on $\sigma$-finite locally compact abelian groups. Let $P$ denote multiplication by
$\mathbf 1_T$ and let $Q=\mathcal F^{-1}\mathbf 1_W\mathcal F$ denote the band-limiting operator, where $T$ and $W$ are subsets of a given group and its dual, respectively.

Smith's uncertainty principle (see Thm.~4.1 in \cite{Smith}) reads as follows: if $f\neq0$ is $(p,\varepsilon)$-supported on $T$ and is $(p,\delta)$-bandlimited to $W$, then
\begin{equation}\label{eq:smith41}
	(\abs T\abs W)^{1/p}\ \ge\ \norm{PQ}_{p\to p}\ \ge\ \frac{1-\ve-\delta}{1+\delta},
\end{equation}
and for $p=2$ the denominator $1+\delta$ is unnecessary.

There are two differences from the setting of \cite{WW} and of Theorem~\ref{main_result}. First, in Smith's
paper the concentration on the Fourier side is measured \emph{in the time domain}, as the
$L^p$-distance to band-limited functions, rather than as the $\ell^p$-tail of the transform. By
Plancherel's theorem the two notions agree for $p=2$, but for $p\neq2$ they need not to be comparable.
Second, as we have already mentioned, Smith's proof relies on translation invariance. However,  Lemma~\ref{lem:smith}
replaces the use of translation invariance by a row-by-row bound. With this modification Smith's argument transfers to special $k$-Hadamard operators. To state a direct analog, we need to define the notion of band-limitedness suitable to our setting.

\begin{definition}\label{def:bandlim}
	Let $A\in\C^{n\times n}$ be special, $\Sigma\subseteq[n]$, $p\in[1,\infty]$ and $\delta\ge0$. A vector
	$v\neq0$ is \emph{$(p,\delta)$-bandlimited to $\Sigma$} if there is $g\in\C^n$ with
	$\supp(Ag)\subseteq\Sigma$ and $\norm{v-g}_p\le\delta\norm v_p$. 
\end{definition}

\begin{prop}\label{prop:smith}
	Let $A\in\C^{n\times n}$ be a special $k$-Hadamard operator, $1\le p\le2$, and $E,\Sigma\subseteq[n]$. Suppose that
	$v \in \C^{n}$ is $(p,\ve)$-supported on $E$ and $(p,\delta)$-bandlimited to $\Sigma$. Then
	\[
	\abs E\,\abs\Sigma\ \ge\ k\Bigl(\frac{1-\ve-\delta}{1+\delta}\Bigr)_+^{p}.
	\]
\end{prop}

\begin{proof}
	Let $X = |E||\Sigma|/k$, as previously. By Lemma~\ref{lem:smith}, $\norm{P_EQ_\Sigma}_{p\to p}\le X^{1/p}$. Suppose that $g$ is as in the Definition ~\ref{def:bandlim}. Since $Q_\Sigma g=g$, we have $P_Eg=P_EQ_\Sigma g$, hence
	\begin{equation}\label{proj}
	v-P_EQ_\Sigma v=(v-P_Ev)+P_E(v-g)+P_EQ_\Sigma(g-v).
	\end{equation}
	This gives 
	\begin{equation}\label{proj2}
	\norm{v-P_EQ_\Sigma v}_p\le(\ve+\delta+X^{1/p}\delta)\norm v_p.
	\end{equation}
	By adding $P_EQ_\Sigma v$ to both sides of \eqref{proj} and combining 
	$\norm{P_EQ_\Sigma v}_p\le X^{1/p}\norm v_p$ with \eqref{proj2}
	we get 
	\[
	1\le\ve+\delta+X^{1/p}(1+\delta),
	\]
	which proves the proposition.
\end{proof}
\subsection{Towards small counterexamples}\label{sec:smallce}

In this subsection we describe results of numerical search aiming at finding smallest counterexample matrices connected to Theorem ~\ref{thm:main}, for $\Phi(k) = k$. The experiment showed existence of such matrices of orders at most 20 (depending on the range of $p$ that we are interested in) while optimization of parameters inside the families from Section ~\ref{counterex} leads to matrices of dimension at least $\sim10^3$. We remark, however, that matrices from Section ~\ref{counterex} are still needed here for the proof of a stronger asymptotic statement in Theorem ~\ref{thm:main}.

For an invertible $A$ with $\abs{a_{ij}}\le 1$,
\[
k_A := \norm{(A^*A)^{-1}}_{\infty\to\infty}^{-1}
= \Bigl(\max_i \sum_j \abs{\bigl((A^*A)^{-1}\bigr)_{ij}}\Bigr)^{-1}
\]
is the largest $k$ for which $A$ is $k$-Hadamard (Definition~\ref{def:kH}).

Let $\ve^{(p)}_E(v)$ and $\eta^{(p)}_\Sigma(Av)$ be as in \eqref{eq:params}, and put
$\ve=\ve^{(p)}_E(v)$, $\eta=\eta^{(p)}_\Sigma(Av)$. Then
$\abs{\supp^p_\ve(v)}\le\abs{E}$ and $\abs{\supp^p_\eta(Av)}\le\abs{\Sigma}$, so any
$A$ and $v$ for which
\begin{equation}\label{eq:violation}
	k_A\bigl(1-\ve^{(p)}_E(v)-\eta^{(p)}_\Sigma(Av)\bigr)^p \ > \ \abs{E}\,\abs{\Sigma}
\end{equation}
violate the conclusion of Theorem~\ref{thm:main} with $\Phi(k)=k$. In all cases below
$E=\Sigma=\{0\}$, so the right-hand side of \eqref{eq:violation} equals $1$.

After normalization, all minimizers with real entries of order
$\le7$ found by the search belong to the two-parameter family
\[
A_m(c,\beta)=\begin{pmatrix}1&\mathbf 1^{\!\top}\\-c\,\mathbf 1&J_m-(1+\beta)I_m\end{pmatrix}\in\R^{(m+1)\times(m+1)},
\]
\[
v=e_0+\frac{c}{m-1-\beta}\sum_{i=1}^me_i ,
\]
with $c,\beta\in[-1,1]$. For $c=\beta=1$ one obtains the $\pm1$ matrix $S_m=A_m(1,1)$ and $v=(m-2,1,\dots,1)$.

The table enclosed below gives every number needed to re-evaluate
the left-hand side (we take $\eta = 0$ everywhere).

\begin{table}[ht]
	\centering\footnotesize
	\setlength{\tabcolsep}{3.5pt}
	\begin{tabular}{lclllll}
		\toprule
		matrix & order & parameters & $k_A$ & $p$ & $\ve_p$ & $k_A(1-\ve_p)^p$\\
		\midrule
		$A_6(c,\beta)$ & 7 & $c=\frac{21943}{10^5}$, $\beta=\frac{44269}{10^5}$ & $\frac{587389}{455731}$ & $1.1$ & $0.20590263$ & $1.000180479$\\
		&&&& $1.5$ & $0.15260172$ & $1.005422517$\\
		&&&& $2$ & $0.11712871$ & $1.004643610$\\
		\midrule
		$A_7(c,\beta)$ & 8 & $c=\frac{4491}{25000}$, $\beta=\frac{10831}{25000}$ & $1.2258893544^{\dagger}$ & $1.1$ & $0.16534485$ & $1.004868096$\\
		&&&& $1.5$ & $0.11499596$ & $1.020631928$\\
		&&&& $2$ & $0.08506921$ & $1.026189952$\\[2pt]
		\midrule
		$S_{13}$ & 14 & $c=\beta=1$ & $\frac{24}{11}$ & $2$ & $0.31147219$ & $1.034335734$\\[2pt]
		\midrule
		$S_{19}$ & 20 & $c=\beta=1$ & $\frac{36}{17}$ & $1.1$ & $0.49079903$ & $1.007933129$\\
		&&&& $1.5$ & $0.35694996$ & $1.091995200$\\
		&&&& $2$ & $0.24837132$ & $1.196355548$\\
		\bottomrule
	\end{tabular}
	\caption{The concrete values of parameters for small counterexamples. $^\dagger$Exact value: $\frac{18684284161352709698}{15241411546905422007}$.}
	\label{tab:cert2}
\end{table}

\begin{remark} 
For matrices $S_m$ all the parameters can be given in closed forms.
With $w=(1,-1,\dots,-1)$ we have $S_m^{\top}S_m=4I+(m-3)ww^{\top}$, and the Sherman--Morrison formula
(\cite{SM}, see also \cite{Hager}, eq.~(2)) gives
\[
(S_m^{\top}S_m)^{-1}=\tfrac14\Bigl(I-\tfrac{m-3}{(m-1)^2}\,ww^{\top}\Bigr),
\]
hence
\[
k_{S_m}=\frac{2(m-1)}{m-2}\quad(m\ge3).
\]
Moreover $S_mv=2(m-1)e_0$ for $v=(m-2,1,\dots,1)$, so $\eta=0$ and
$\ve_p^p=\dfrac{m}{(m-2)^p+m}$. 

Within this family, the number $m=13$ is the smallest for which we obtained a contradiction for $p=2$. For $m=19$ the violation of the uncertainty principle was checked for $p$ arbitrarily close to one within $50$-digit arithmetic.
\end{remark}

\section*{Acknowledgements}
The article was written in collaboration with FWF Project P 34414, \textit{Radial Variation, Green's Mappings and Brownian Motion}.

\section*{Declaration of the use of generative AI}

The author used Claude (Anthropic) to assist with the development of mathematical arguments and the preparation of computational experiments. Claude found counterexamples from Sections 3 and 5, prepared the baseline for numerical experiments mentioned in Section 6 and assisted in the proof of Theorem 1.8. 

All mathematical statements and proofs were independently reviewed and verified by the author, who takes full responsibility for the contents of this paper.

\bibliography{biblio} 
\bibliographystyle{alpha}

\end{document}